\documentclass[a4paper,12pt]{amsart}

\usepackage{amssymb, amscd}
\usepackage{latexsym,epsfig}
\usepackage[all]{xy}
\usepackage{pst-all}
\numberwithin{equation}{section}
\def\today{\ifcase\month\or Jan\or Febr\or  Mar\or  Apr\or May\or Jun\or  Jul\or Aug\or  Sep\or  Oct\or Nov\or  Dec\or\fi \space\number\day, \number\year}

\newcommand{\AAA}{\mathbb A}
\newcommand{\CC}{\mathbb C}

\newcommand{\FF}{\mathbb F}

\newcommand{\PP}{\mathbb P}
\newcommand{\QQ}{\mathbb Q}

\newcommand{\VV}{\mathbb V}
\newcommand{\WW}{\mathbb W}
\newcommand{\ZZ}{\mathbb Z}

\numberwithin{equation}{section}
\newcommand{\Vladut}{Vl\u adu\c t}

\newtheorem{theorem}{Theorem}[section]

\newtheorem{conjecture}[theorem]{Conjecture}

\newtheorem{definition-lemma}[theorem]{Definition-Lemma}
\theoremstyle{definition}

\newtheorem{example}[theorem]{Example}

\theoremstyle{remark}

\begin{document}

\title[]{Counting Curves over Finite Fields}
\author{Gerard van der Geer}
\address{Korteweg-de Vries Instituut, Universiteit van
Amsterdam, Postbus 94248,
1090 GE  Amsterdam, The Netherlands.}
\email{G.B.M.vanderGeer@uva.nl}

\begin{abstract}
This is a survey on recent results on counting of curves over finite fields.
It reviews various results on the maximum number of points on a curve
of genus $g$ over a finite field of cardinality $q$, but the main
emphasis is on results on the Euler characteristic of the 
cohomology of local systems on moduli
spaces of curves of low genus and its implications for modular forms.
\end{abstract}

\maketitle
\begin{section}{Introduction}\label{sec-intro}
Reduction modulo a prime became a standard method for studying equations
in integers after Gauss published his Disquisitiones Arithmeticae in 1801.
In \S 358 of the Disquisitiones Gauss counts the number of 
solutions of the cubic Fermat equation $x^3+y^3+z^3=0$ modulo a prime $p$
and finds for a prime $p \not\equiv 1 (\bmod 3)$ always $p+1$ points on the
projective curve, while
for a prime $p\equiv 1 (\bmod 3)$ 
the number of points
equals $p+1+a$, if one writes $4p=a^2+27 b^2$ with $a\equiv 1 (\bmod 3)$ 
and he notes that $|a|\leq 2 \sqrt{p}$.
But although Galois introduced finite fields in 1830 and algebraic curves
were one of the main notions in 19th century mathematics, one had to wait
till the beginning of the 20th century before algebraic curves over finite
fields became an important topic for mathematical investigation. 
Artin considered
in his 1924 thesis (already submitted to Mathematische Zeitschrift in 1921)
the function fields of hyperelliptic curves defined over a finite field
and considered for such fields a zeta function $Z(s)$ that is an analogue of 
the Riemann zeta function and of the Dedekind zeta function for number fields.
He derived a functional equation
for them and formulated an analogue of the Riemann hypothesis that says
that the zeros of the function of $t$ obtained by substituting 
$t=q^{-s}$ in $Z(s)$ have absolute value $q^{-1/2}$.
In 1931 Friedrich Karl Schmidt brought a more geometric approach 
by writing the zeta function
for a smooth absolutely irreducible 
projective curve $C$ over a finite field ${\FF}_q$ as the 
generating function for the number of rational points 
$c(n)=\# C({\FF}_{q^n})$
over extension fields as
$$
Z(t)= {\rm exp} (\sum_{n=1}^{\infty} c(n) \frac{t^n}{n})\, ,
$$
which turns out to be a rational function of $t$ of the form
$$
Z(t)=\frac{P(t)}{(1-t)(1-qt)}
$$
for some polynomial $P\in {\ZZ}[t]$ 
of degree $2g$ with $g$ the genus of the curve.
He observed that the functional equation $Z(1/qt)=q^{1-g}t^{2-2g}Z(t)$
is a consequence of the theorem of Riemann-Roch. 
A couple of years later (1934)  Hasse proved the Riemann hypothesis 
for elliptic curves over finite fields using correspondences. 
The proof appeared in 1936, see \cite{Hasse}.
Deuring observed
then that to extend this result to curves of higher genus one needed
a theory of algebraic correspondences over fields of arbitrary 
characteristic. This was at the time that the need was felt 
to build algebraic geometry on a more solid base that would allow one to
do algebraic geometry over arbitrary fields. 
Weil was one of those who actively pursued this goal.
Besides doing foundational work, he also exploited the analogy 
between geometry in characteristic zero
and positive characteristic by extending an inequality on correspondences of
Castelnuovo and Severi to positive characteristic and deduced around 1940 
the celebrated Hasse-Weil inequality
$$
| \#C({\FF}_q)-(q+1)| \leq 2g \, \sqrt{q}
$$
for the number of rational points on a smooth absolutely irreducible 
projective curve $C$ of genus $g$ over a finite field ${\FF}_q$ (\cite{Weil}).

Geometry entered the topic more definitely 
when Weil applied the analogy with
the Lefschetz fixed point theorem, which expresses the number of fixed points
of a map on a compact manifold in terms of the trace of the induced
map on (co-)homology,  to the case of the Frobenius morphism on a projective
variety over a finite field and formulated in 1949 
the famous `Weil Conjectures'
on zeta functions of varieties over finite fields. Dwork set the first step
by proving  the
rationality of the zeta function in 1960.

Grothendieck's revolution in algebraic geometry in the late 1950s 
started a new era in which it was possible to do algebraic geometry 
on varieties over  finite fields. It also led to the construction
of \'etale cohomology, which made it possible to carry out the analogy
envisioned by Weil.
The first milestone in this new era was Deligne's 
completion of the proof of the Weil conjectures in 1974. 

Among all these developments the theme of curves over finite fields  
was pushed to the background, though there was progress.
In 1969 Stepanov showed a new approach in \cite{Stepanov}
to deriving the Hasse-Weil bound for hyperelliptic curves
by just using Riemann-Roch; Stark used it to get a somewhat stronger
bound than Hasse-Weil for hyperelliptic curves over a prime
field ${\FF}_p$, see \cite{Stark}.
Stepanov's method was elegantly extended by Bombieri in \cite{Bombieri} 
to prove the Hasse-Weil bound in the general case.

The return of curves over finite fields to the foreground around 1980
was triggered by an outside impulse,
namely from coding theory. Goppa observed that one could construct
good codes by evaluating meromorphic functions on a subset of 
the points of the projective line, where ``good" meant 
that they reached the so-called Gilbert-Varshamov bound. 
He then realized that 
this could be generalized by evaluating meromorphic functions 
on a subset of the rational points of a higher genus curve,
that is, by associating a code to a linear system on a curve
over a finite field, cf.\ \cite{Goppa}. 
The quality of the code depended on the number of rational points
of the curve. In this way it drew attention to the question
how many rational points a curve of given genus $g$ over a finite
field ${\FF}_q$ of given cardinality $q$ could have.
The 1981 paper by Manin \cite{Manin1} explicitly asks in the title 
for the maximum number of points on a curve over ${\FF}_2$.
Thus the question emerged how good the Hasse-Weil bound was. 
\end{section}
\begin{section}{The Maximum Number of Points on a Curve over a Finite Field}
In \cite{Ihara} Ihara employed a simple idea to obtain a better estimate 
than the Hasse-Weil bound for
the number of rational points on a curve over a finite field ${\FF}_q$. 
The idea is to write 
$$
\# C({\FF}_{q^n})= q^n+1-\sum_{i=1}^{2g} \alpha_i^n\, ,
$$
with $\alpha_i$ the eigenvalues of Frobenius on $H^1_{\rm et}(C,{\QQ}_{\ell})$
for $\ell$ different from the characteristic, 
and to note that $\#C({\FF}_{q})  \leq \# C({\FF}_{q^2})$.
By using the Cauchy-Schwartz inequality for the $\alpha_i$
he found the improvement
$$
\# C({\FF}_q)\leq q+1+[(\sqrt{(8q+1)g^2+4(q^2-q)g} \, - g)/2]\, .
$$
For $g> (q-\sqrt{q})/2$ this is better than the Hasse-Weil bound.
Instead of just playing off $\#C({\FF}_q)$ against $\#C({\FF}_{q^2})$,
one can use the extensions of ${\FF}_q$ of all degrees, 
and a systematic analysis (in \cite{D-V})
due to Drinfel'd  and \Vladut\  led to an asymptotic upper bound
for the quantity
$$
A(q):= \lim \sup_{g \to \infty} N_q(g)/g \, ,
$$
which was introduced by Ihara, with $N_q(g)$ as usual defined as
$$
N_q(g):= \max \{ \# C({\FF}_q): g(C)=g \}\, ,
$$
the maximum number of rational points on a smooth absolutely irreducible 
projective curve of genus $g$ over ${\FF}_q$. 
The resulting asymptotic bound is
$$
A(q)\leq \sqrt{q} -1.
$$
As we shall see below, this is sharp for $q$ a square.

The systematic study of $N_q(g)$ was started by Serre in the 1980s. 
He showed in \cite{Serre1} that by using some arithmetic 
the Hasse-Weil bound can be improved slightly
to give
$$
| \# C({\FF}_q)-(q+1)| \leq g \lfloor  2\sqrt{q} \rfloor\, ,
$$
as opposed to just $\leq \lfloor 2g\sqrt{q} \rfloor$. Serre applied 
the method of `formules explicites' from number theory 
to the zeta functions of curves over finite fields to get better upper bounds. 
An even trigonometric polynomial
$$
f(\theta)= 1+2 \sum_{n\geq 1} u_n \cos n\theta
$$
with real coefficients $u_n \geq 0$ such that $f(\theta)\geq 0$
for all real $\theta$ gives an estimate for $\#C ({\FF}_q)$ of the form
$$
\# C({\FF}_q)\leq a_f g + b_f,
$$
with $g$ the genus of $C$ and $a_f$ and $b_f$ defined by setting
$\psi=\sum_{n\geq 1} u_nt^n$ and
$$
a_f= \frac{1}{\psi(1/\sqrt{q})} \quad {\rm and} \quad
b_f=1+ \frac{\psi({\sqrt{q})}}{\psi(1/\sqrt{q})}.
$$
Oesterl\'e found the solution to the problem of finding the optimal choices for
the functions~$f$, see \cite{Schoof}. For $g> (q-\sqrt{q})/2$
these bounds are better than the Hasse-Weil bound.

Curves that reach the Hasse-Weil upper bound are called {\sl maximal curves}.
In such a case $q$ is a square and $g \leq (q-\sqrt{q})/2$. 
Stichtenoth and Xing conjectured 
that for maximal curves over ${\FF}_q$ 
one either has that $g=(q-\sqrt{q})/2$
or $g \leq (\sqrt{q}-1)^2/4$, and after they made considerable progress
towards it, see \cite{SX}, the conjecture was proved by Fuhrmann and Torres,
cf.\ \cite{FT}. 
In this direction it is worth mentioning a recent result of Elkies, Howe and 
Ritzenthaler (\cite{EHR}) that gives a bound on the genus of curves whose Jacobian 
has Frobenius eigenvalues in a given finite set.

Stichtenoth conjectured that
all maximal curves over ${\FF}_{q^2}$ 
are dominated by a `hermitian' curve defined by an
equation
$$
x^{q+1}+y^{q+1}+z^{q+1}=0\, .
$$
This curve is of genus $q(q-1)/2$ and has $q^3+1$ 
rational points over ${\FF}_{q^2}$. This conjecture was disproved 
by Giulietti and Korchm\'aros, who exhibited a counterexample over ${\FF}_{q^6}$,
see \cite{G-K}. 
In \cite{RS} R\"uck and Stichtenoth proved that maximal curves with $g=(q-\sqrt{q})/2$ 
are isomorphic to the hermitian curve. 

For a curve $C$ over a finite field ${\FF}_q$ the quantity
$$
\delta:= (q+1+g[2\sqrt{q}])-\# C({\FF}_q)
$$
is called the {\sl defect}.
The result of Fuhrmann and Torres proves the non-existence of curves
with a small defect. Many more results excluding curves with small
defects have been obtained
by various arithmetic and geometric methods, 
see especially work of Howe and Lauter; we refer to 
the papers \cite{HL1,HL2,HL3, L1,L2,Serre4}. 

However, testing how good the
resulting upper bounds on $N_q(g)$ are,
can only be done by providing a curve with a 
number of points that reaches or comes close to this upper bound; that is,
by constructing a curve with many points.

In \cite{Serre1, Serre2, Serre3} 
Serre listed the value of $N_q(g)$ for small values of 
$q$ and $g$ or a small interval in which $N_q(g)$ lies when the value of
$N_q(g)$ was not known. Wirtz extended in \cite{Wi} these tables for small 
$q$ that are powers of $2$ and $3$ by carrying out a computer search
in certain families. (His table is reproduced in \cite{vdG-vdV2}, p.\ 185.) 
In the 1990s the challenge to find curves over finite fields with many
points, that is, close to the best upper bound for $N_q(g)$, attracted
a lot of interest. 
In 1996 van der Geer and van der Vlugt published
`Tables for the function $N_q(g)$' that listed intervals for the 
function $N_q(g)$ for $1\leq g \leq 50$ and $q$ a small power of $2$ or $3$. 
These tables were regularly updated and published on a website. In 1998
the tables were replaced by a new series of tables 
(`Tables of Curves with Many Points'),
one of the first of which was published in  
Mathematics of Computation \cite{vdG-vdV3}, 
and it was regular updated on a website. 
In a series of papers (see \cite{NX1} and the references there)
Niederreiter and Xing efficiently 
applied methods from class field theory to construct curves with many points,
resulting in many good entries in the tables. 
Other methods, like the ones used in \cite{vdG-vdV2, vdG2},
employed fibre products of Artin-Schreier curves
or were based on coding theory, see \cite{vdG-vdV3} and the references
given there. The resulting tables were the joint effort of many people.
The last update was dated October 2009; after that 
the tables were replaced by a new website, {\tt www.manypoints.org},
an initiative of van der Geer, Howe, Lauter and Ritzenthaler, 
where new records can be registered. 
At the end of this review we include a copy
of a recent version of the tables for small powers of $2$ and~$3$
and $1\leq g \leq 50$.
As the reader will see, the intervals for $N_q(g)$ are still 
quite large for many pairs $(g,q)$.

As mentioned above, the result of Drinfeld and \Vladut\ 
led to the asymptotic bound $A(q)\leq \sqrt{q}-1$. 
For $q$ a square, Ihara and independently Tsfasman, \Vladut\  and Zink showed in 
\cite{TVZ,Ihara} that modular curves have many rational points 
and that one can use this to prove
$$
A(q) \geq \sqrt{q}-1\, ,
$$
so that $A(q)=\sqrt{q}-1$ for $q$ a square. It came as a surprise 
in 1995 when Garcia and Stichtenoth  came 
forward (see \cite{GS1})
with a tower over ${\FF}_{q^2}$ ($q$ an arbitrary prime power)
$$ 
\ldots C_i \to C_{i-1} \to \cdots \to C_2 \to C_1
$$ 
of Artin-Schreier curves defined over ${\FF}_{q^2}$ 
by a simple recursion 
with 
$$
\lim_{i\to \infty} g(C_i)=\infty \qquad {\rm and} \qquad 
\lim_{i\to \infty} \frac{\# C_i({\FF}_{q^2})}{g(C_i)} =q-1\, .
$$
The simple recursion starts with ${\PP}^1$ over ${\FF}_{q^2}$ with function field $F_1={\FF}_{q^2}(x_1)$ and defines Artin-Schreier extensions 
$F_n$ by with $F_{n+1}=F_n(y_{n+1})$ given by
$y_{n+1}^q+y_{n+1}= x_n^{q+1}$ with $x_{n+1}:=y_{n+1}/x_n$ for $n\geq 1$.  
This has stimulated much research. 
Elkies has shown in \cite{Elkies} 
that this tower is in fact a tower of modular curves.
Over fields the cardinality of which is not a square it is more
difficult to find good towers. There are towers resulting from
class field theory, see for example \cite{NX1}. In a paper from 1985,
\cite{Zink},
Zink used certain degenerate 
Shimura surfaces to construct a tower over ${\FF}_{p^3}$
for $p$ prime  with limit  $\#(C_i({\FF}_{p^3}))/g(C_i) \geq
2(p^2-1)/(p+2)$.
The first good explicit wild tower in the non-square case
was a tower of Artin-Schreier covers
over ${\FF}_8$ with limit $3/2$, see \cite{vdG-vdVt}.
This has been generalized by Bezerra, Garcia and Stichtenoth to towers over
${\FF}_{q}$ with $q$ a cube. If $q=\ell^3$ with $\ell$  prime power they
deduce that
$$
A(\ell^3)\geq \frac{2(\ell^2-1)}{\ell+2}\, ,
$$
and this was extended again in \cite{BBGS} to all nonprime finite fields ${\FF}_q$.
For a detailed review of the progress on towers we refer to the 
paper by Garcia and Stichtenoth \cite{GS2}. 

The question of the maximum number of points on a curve of given genus $g$ 
over a finite field ${\FF}_q$ is just one small part of the question which
values the number of points on a curve of genus $g$ over ${\FF}_q$ can have.
The answer can take various forms. One answer is in \cite{AS},
where it is shown that for sufficiently large genus, every 
value in a small interval $[0, c]$ is assumed.
But, more precisely, one may ask which values are assumed and how often if the
curve varies through the moduli space of curves of genus $g$ defined 
over ${\FF}_q$. 
The answer could be presented as a list of all possible Weil polynomials
with the frequencies with which they occur. 
The question arises how to 
process all the information contained in such a list. 
For example, take the case of elliptic curves over a 
finite prime field ${\FF}_p$.
Each isomorphism class $[E]$ of elliptic curves defined over ${\FF}_p$ defines a
pair $\{ \alpha_E, \bar{\alpha}_E\}$ of algebraic integers with
$\# E({\FF}_p)=p+1-\alpha_E-\bar{\alpha}_E$. We can study the 
weighted `moments' 
$$
\sigma_k(p):=
-\sum \frac{ \alpha_E^k+\alpha_E^{k-1}\bar{\alpha}_E+\cdots + \bar{\alpha}_E^k}{
\# {\rm Aut}_{{\FF}_p}(E)}\, ,
$$
where the sum is over a complete set of
representatives $E$ of all the isomorphism classes of elliptic
curves over ${\FF}_p$. 
For odd $k$ the answer is $0$ due to the fact that the contribution of 
an elliptic curve and its $-1$-twist  cancel. For $k=0$ we find $-q$,
while for even $k$ with $2\leq k \leq 8$ we find $1$. But for $k=10$
we find for $p=2,3,5,7,11$ the following values
\smallskip
\vbox{
\bigskip\centerline{\def\quad{\hskip 0.6em\relax}
\def\quod{\hskip 0.5em\relax }
\vbox{\offinterlineskip
\hrule
\halign{&\vrule#&\strut\quod\hfil#\quad\cr
height2pt&\omit&&\omit&&\omit&&\omit&&\omit&&\omit&\cr
& $p$ && $2$ && $3$ && $5$ && $7$ && $11$ &\cr
\noalign{\hrule}
& $\sigma_{10}$ && $-23$ && $253$ && $4831$ && $-16743$ && $534613$ &\cr
} \hrule}
}}

\bigskip
Many readers will not fail to notice that the numbers appearing 
here equal $\tau(p)+1$
with $\tau(p)$ the $p$th Fourier coefficient of the celebrated modular form 
$\Delta=\sum \tau(n)q^n$ of weight $12$ on ${\rm SL}(2,{\ZZ})$ 
with Fourier development
$$
\Delta= q-24\, q^2+252\, q^3-1472 \, q^4+4830 \, q^5-6048 \, q^6-16744 \, q^7 +\cdots
$$
where the reader will hopefully forgive us for having used the customary
$q=e^{2\pi i \tau}$.
This hints at treasures hidden in such frequency lists of Weil polynomials
and the rest of this survey paper is dedicated to this phenomenon.
(For a different view on such statistics we refer to \cite{BG}.)

\end{section}
\begin{section}{Varieties over the Integers}
A customary approach for studying algebraic varieties begins by trying to 
 calculate their cohomology.
For a variety defined over a finite field ${\FF}_q$ we can extract
a lot of information on the cohomology 
by counting rational points of the variety
over the extension fields ${\FF}_{q^r}$. 
The connection is through the Lefschetz trace formula
which says that the number of points equals the trace of
the Frobenius morphism on the rational Euler characteristic
of the variety.
And for a variety defined
over the integers we can look at its reduction modulo a prime and
then count rational points over extension fields ${\FF}_{p^r}$. 
This characteristic $p$ information
can then be pieced together to find cohomological
information about the variety in characteristic zero,
more precisely, about the cohomology as a representation
of the absolute Galois group of the rational numbers.

Let us look at proper varieties defined over the integers with 
good reduction everywhere. 
The first examples are given by projective space and Grassmann varieties.
For projective space we have $\#{\PP}^n({\FF}_q)=q^n+q^{n-1}+\cdots+1$
and for the Grassmann variety  $G(d,n)$ 
of $d$-dimensional projective linear subspaces of ${\PP}^n$ we have
$$
\# G(d,n)=  \left[ {n+1\over d+1}\right]_q :=
\frac{(q^{n+1}-1)(q^{n+1}-q) \cdots (q^{n+1}-q^d)}{(q^{d+1}-1)(q^{d+1}-q)\cdots (q^{d+1}-q^d)}\, .
$$
In fact, for these varieties we have a cell decomposition and we 
know the class in the Grothendieck group of varieties. 
Recall that if $k$ is a perfect field $k$
and ${\rm Var}_k$ is the category of algebraic
varieties over $k$, then the Grothendieck group $K_0({\rm Var}_k)$
of varieties over $k$ is, by definition, the free abelian group 
generated by the symbols $[X]$ with $X$ an object of ${\rm Var}_k$
modulo the two relations
i) $[X]=[Y]$ whenever $X \cong Y$;
ii) for every closed subvariety $Z$ of $X$ we have
$[X]=[Z]+[X\backslash Z]$.
The class of the affine line ${\AAA}^1$ is denoted by ${\mathbf L}$ and called the
Lefschetz class. For example, 
for the projective space ${\PP}^n$ and the Grassman variety $G(d,n)$ 
we find
$$
[{\PP}^n]= {\mathbf L}^n+{\mathbf L}^{n-1}+\cdots + 1
\qquad \text{and} \qquad 
[G(d,n)]=\left[ {n+1\over d+1}\right]_{\mathbf L} \, .
$$

For varieties $X$ like  projective spaces and Grassmannian varieties,
where we know a cell decomposition,
we find that there exists 
a polynomial $P\in {\ZZ}[x]$ such that $\# X({\FF}_q)=P(q)$
for every finite field ${\FF}_q$.
Conversely, one can ask how much we can learn about  a proper
smooth variety defined over the integers
by counting the number of ${\FF}_q$-rational points for many fields ${\FF}_q$.

For example, if we find that there exists such a polynomial $P$ with
$\# X({\FF}_q)=P(q)$, what do we know? 
There is a theorem by 
van den Bogaart and Edixhoven (\cite{vdB-E}) 
which says that for a proper variety we  then know the 
$\ell$-adic \'etale cohomology for all $\ell$ of $X_{\overline{\QQ}}$ as a
representation of the absolute Galois group  of the rational 
numbers: the cohomology is a
direct sum of copies of the cyclotomic representation 
${\QQ}_{\ell}(-i)$ 
in degree $2i$ and zero
in odd degrees; morover the number of copies of ${\QQ}_{\ell}(-i)$
is given by the $i$th coefficient of $P$. (For non-proper spaces then the result holds
for the Euler characteristic in a suitable Grothendieck
group.) Note that the realization of
the motive ${\mathbf L}^i$ as a Galois representation equals ${\QQ}_{\ell}(-i)$.

The spaces ${\PP}^n$ and the Grassmann varieties 
are moduli spaces as they parametrize linear subspaces of projective
space. The first further examples of varieties defined over the integers
with everywhere good reduction 
are also moduli spaces, the moduli spaces ${\mathcal M}_g$ of curves of
genus $g$ and the moduli spaces 
${\mathcal A}_g$ of principally polarized abelian varieties
of dimension $g$. More generally, there are the moduli
spaces ${\mathcal M}_{g,n}$
of $n$-pointed curves of genus $g$ and their Deligne-Mumford
compactifications $\overline{\mathcal M}_{g,n}$ of stable $n$-pointed curves. 
All these spaces, 
${\mathcal A}_g$, ${\mathcal M}_g$, ${\mathcal M}_{g,n}$ 
and $\overline{\mathcal M}_{g,n}$, are
Deligne-Mumford stacks defined over the integers and smooth over ${\ZZ}$.
The spaces $\overline{\mathcal M}_{g,n}$ 
are also proper over ${\ZZ}$. 
These spaces constitute the most intriguing 
series of varieties (or rather Deligne-Mumford stacks)
over the integers with everywhere good reduction. 
In the last two decades our knowledge about them has increased dramatically,
but clearly so much remains to be discovered.

While the cohomology of projective space and the Grassmann varieties
is a polynomial in ${\mathbf L}$ (or ${\QQ}(-1)$),
it is unreasonable to expect the same for the moduli space ${\mathcal A}_g$
and ${\mathcal M}_{g,n}$.
In fact, we know that over the complex numbers ${\mathcal A}_g$
can be described as a quotient 
${\rm Sp}(2g,{\ZZ})\backslash \mathfrak{H}_g$, 
with $\mathfrak{H}_g$ the Siegel upper half space (see below), and
that modular forms are supposed to contribute to its cohomology.
In fact, 
the compactly supported cohomology possesses a mixed Hodge structure and
cusp forms of weight $g+1$ (see next section) contribute to the
first step in the Hodge filtration on middle-dimensional cohomology. 
Since we know that for large $g$
there exist non-trivial cusp forms of this weight (e.g.\ $g=11$) 
this shows that the cohomology is not
so simple. In fact, we know that the cohomology can be described in
terms of automorphic forms on the symplectic group ${\rm Sp}(2g)$. 
Despite this, for low values of $g$ and $n$ the cohomology of
${\mathcal M}_{g,n}$ can be
a polynomial in ${\mathbf L}$. For example, for $g=1$ and $n\leq 9$,
for $g=2$ with $n \leq 7$ and $g=3$ for $n\leq 7$ we have explicit 
polynomial formulas for the number of points over finite fields,
and hence for the Euler characteristic of the
moduli space ${\mathcal M}_{g,n}\otimes {\FF}_q$ as a polynomial in ${\mathbf L}$,
see Getzler \cite{Getzler99} and Bergstr\"om \cite{Bergstroem}.

If one does not find an explicit polynomial in $q$ that gives
the number of ${\FF}_q$-rational points on our moduli space over ${\FF}_q$,
one nevertheless might try to count the number of ${\FF}_q$-rational points 
of ${\mathcal A}_g\otimes {\FF}_p$ to get information on the Euler characteristic 
of the cohomology.
Since ${\mathcal A}_g \otimes {\FF}_p$ (or ${\mathcal M}_g\otimes {\FF}_p$)
is a moduli space its points are
represented by objects (abelian varieties or curves) 
and the first question then is 
how to represent the objects parametrized by ${\mathcal
A}_g$ (or ${\mathcal M}_g$). 
For $g=1$ this is clear. If we make a list of all elliptic
curves defined over ${\FF}_q$ up to isomorphism over ${\FF}_q$
and calculate for each such elliptic curve the number of ${\FF}_q$-rational
points we should be able to calculate $\# {\mathcal M}_{1,n}({\FF}_q)$
for all $n\geq 1$. (This has to be taken with a grain of 
salt as ${\mathcal M}_{1,n}$ is a stack and not a variety;
this aspect is taken care of by taking into account the automorphism
groups of the objects.)
This is the approach we shall take in the next section.
\end{section}
\begin{section}{Counting Points on Elliptic Curves}
Hasse proved in 1934 that the number of rational points on an elliptic
curve $E$ defined over a finite field ${\FF}_q$ can be given as
$$
\# E({\FF}_q)= q+1-\alpha-\bar{\alpha}
$$
with $\alpha=\alpha_E$ an algebraic integer with $\alpha \bar{\alpha}=q$. 
Isomorphism classes of elliptic curves over the algebraic closure 
$\overline{\FF}_q$ are given by
their $j$-invariant $j(E)$; over the field ${\FF}_q$ this is no longer true
due to automorphisms of the curve. But for any given value of 
$j \in {\FF}_q$ there is an elliptic curve $E_j$ defined over 
${\FF}_q$ and the ${\FF}_q$-isomorphism classes of elliptic curves 
defined over ${\FF}_q$ with this $j$-invariant
correspond $1-1$ with the elements of the pointed set 
$H^1({\rm Gal}_{\overline{\FF}_q/{{\FF}_q}}, {\rm Isom}(E_j))$
with ${\rm Isom}$ the group of ${\FF}_q$-automorphisms of the genus $1$ curve
underlying $E$. For each $E_j$ this set contains at least two elements.
Nevertheless, a given value of $j\in {\FF}_q$ contributes just
$$
\sum_{E/{\FF}_q / \cong_{{\FF}_q},j(E)=j}   
\frac{1}{\# {\rm Aut}_{{\FF}_q}(E)}=1
$$
to the number of elliptic curves defined over ${\FF}_q$ up to 
${\FF}_q$-isomorphism, if we count them in the right way, that is, with
weight $1/\# {\rm Aut}_{{\FF}_q}(E)$, see \cite{vdG-vdV1} for a proof.
 
We are interested in how the $\alpha$ vary over the whole $j$-line. 
To this end one considers the moments of the $\alpha_E$
$$
\sigma_a(q):=
-\sum_{E/{\FF}_q /\cong_{{\FF}_q}} \frac{\alpha_E^a+\alpha_E^{a-1}\bar{\alpha}_E+\cdots + \bar{\alpha}_E^{a}}{\# {\rm Aut}_{{\FF}_q}(E)}\, ,
\eqno(1)
$$
where the sum is over all elliptic curves defined over ${\FF}_q$
up to isomorphism over ${\FF}_q$ and $a$ is a non-negative integer. 
For odd $a$ one finds zero, due to the fact that the contributions 
of a curve and its $-1$-twist cancel. But for even $a>0$ one finds something 
surprising and very interesting at which we hinted at the end of Section 2: 
for a prime $p$ we have
$$
\sigma_a(p)= 1+\text{Trace of $T(p)$ on $S_{a+2}({\rm SL}(2,{\ZZ}))$}
\eqno(2)$$
with $S_k({\rm SL}(2,{\ZZ}))$ the space of cusp forms of weight $k$
on ${\rm SL}(2,{\ZZ})$ and $T(p)$ the Hecke operator associated to $p$ 
on this space. 
Recall that a modular form of weight $k$ on ${\rm SL}(2,{\ZZ})$ is a
holomorphic function $f: \mathfrak{H} \to {\CC}$ on the upper half plane
$\mathfrak{H}=\{ \tau \in {\CC}: {\rm Im}(\tau) >0 \}$ 
of ${\CC}$ that satisfies
$$
f(
\frac{a\tau+b}{c\tau+d})
=(c\tau+d)^k f(\tau) \qquad
\text{for all $\tau \in \mathfrak{H}$ 
and 
$\left( \begin{matrix} a & b \\ c & d \\ \end{matrix} \right) 
\in {\rm SL}(2, {\ZZ})$ }\, ,
$$
in particular, it satisfies $f(\tau+1)=f(\tau)$ 
and thus admits a Fourier development
$$
f=\sum_{n} a(n)\, e^{2\pi i n \tau},
$$
and we require that $f$ be holomorphic at infinity, i.e. $a(n)=0$ for $n<0$.
A {\sl cusp form} is a modular form with vanishing constant term $a(0)=0$.
The modular forms of given weight $k$ form a vector space 
$M_k({\rm SL}(2,{\ZZ}))$ of finite dimension; this dimension
is zero for $k$ negative or odd and equals 
$[k/12]+1$ for even $k\not\equiv 2 (\bmod 12)$ and $[k/12]$ for even 
$k\equiv 2 (\bmod 12)$. The subspace $S_k({\rm SL}(2,{\ZZ}))$ of cusp forms
of weight $k$ is of codimension $1$ in $M_k({\rm SL}(2,{\ZZ}))$
if the latter is nonzero.
One has an algebra of Hecke operators $T(n)$ with $n\in {\ZZ}_{\geq 1}$
operating on $M_k({\rm SL}(2,{\ZZ}))$ and $S_k({\rm SL}(2,{\ZZ}))$
and there is a basis of common eigenvectors, called eigenforms,
for all $T(n)$ with the property that $T(n)f= a(n)f$ for such an eigenform
$f$ if one normalizes these such that $a(1)=1$.
Modular forms belong to the most important objects in arithmetic
algebraic geometry and number theory. 
It may come as a surprise that we can obtain information 
about modular forms, that are holomorphic functions on $\frak{H}$, 
by counting points on elliptic curves over finite fields.

On the other hand, our knowledge of modular forms on ${\rm SL}(2,{\ZZ})$ 
is extensive. Since a product of modular forms of weight $k_1$ and $k_2$ is a
modular form of weight $k_1+k_2$ one obtains a graded algebra
$\oplus_k M_k({\rm SL}(2,{\ZZ}))$ of modular forms 
on ${\rm SL}(2,{\ZZ})$ and it is the polynomial
algebra generated by
the Eisenstein series $E_4$ and $E_6$, explicit modular forms 
of weight $4$ and $6$.
This does not tell us much about the action of the Hecke operators, but
the fact is that
one has a  closed formula for the trace of the Hecke operator $T(n)$ on the 
space $S_k({\rm SL}(2,{\ZZ}))$ for even $k>0$:
$$
{\rm Trace}(T(n))= 
-\frac{1}{2} \sum_{t=-\infty}^{\infty} P_k(t,n) \, H(4n-t^2) -
\frac{1}{2} \sum_{dd'=n} \min(d,d')^{k-1}\, ,
$$
where $P_k(t,n)$ is the coefficient of $x^{k-2}$ in the Taylor series 
of $(1-tx+nx^2)^{-1}$ and $H(n)$ is a class number defined as follows. 
For $n<0$ we put $H(n)=0$; furthermore $H(0)=-1/12$, 
while for $n>0$ we let $H(n)$ be the number of 
${\rm SL}(2,{\ZZ})$-equivalence classes
of positive definite binary quadratic forms 
$ax^2+bxy+cy^2$ of discriminant $b^2-4ac=-n$ with the forms 
equivalent to $x^2+y^2$ (resp.\ to $x^2+xy+y^2$) 
counted with weight $1/2$ (resp.\ $1/3$).
So in view of all we know, 
we do not gain new information from our counts of points over finite fields.
\bigskip

The fact that we can obtain information on modular forms 
by counting points over finite fields
illustrates two ideas: the idea of Weil that counting points on varieties 
over finite fields gives the trace of Frobenius on the  cohomology
and the idea that the cohomology of the variety obtained by reducing a variety 
defined over the integers modulo $p$ reflects aspects of the
cohomology of the variety over the integers. Moreover, it shows that 
modular forms are cohomological invariants.

Note that the expression $\sigma_a(p)$ in (1) 
is the sum over all elliptic curves 
defined over ${\FF}_p$ up to isomorphism 
of the negative of the trace 
of Frobenius on the $a$th symmetric power of the cohomology
$H^1_{\rm et}(E\otimes \overline{\FF}_p, {\QQ}_{\ell})$ with $\ell$ a prime different from $p$.

One thus is led to look at the local system 
${\VV}:=R^1\pi_* {\QQ}_{\ell}$ on
the moduli space ${\mathcal A}_1$ of elliptic curves with 
$\pi: {\mathcal X}_1 \to {\mathcal A}_1$ the universal family of 
elliptic curves. This is a local system of vector spaces
with fibre over $[E]$ equal to $H^1(E,{\QQ}_{\ell})$.
Note that $\pi: {\mathcal X}_1 \to {\mathcal A}_1$ 
is defined over ${\ZZ}$. 

For even $a>0$ we look at the local system ${\VV}_a=
{\rm Sym}^a({\VV})$; this is a local system of rank $a+1$ 
on ${\mathcal A}_1$ with fibre
${\rm Sym}^a(H^1(E,{\QQ}_{\ell}))$ over a point  $[E]$ of
 the base ${\mathcal A}_1$. 
It is in the cohomology of ${\VV}_a\otimes {\CC}$ 
over ${\mathcal A}_1 \otimes {\CC}$ 
that we find the modular forms. In fact, 
a famous theorem of Eichler and Shimura says that
$$
H^1({\mathcal A}_1 \otimes {\CC}, {\VV}_a \otimes {\CC})\cong 
S_{a+2}({\rm SL}(2,{\ZZ}))\oplus
\overline{S}_{a+2}({\rm SL}(2,{\ZZ}))\oplus {\CC}
\eqno(3)
$$
So the space $S_{a+2}({\rm SL}(2,{\ZZ}))$ of cusp forms of weight $a+2$ 
and its complex conjugate constitute
this cohomology, except for the summand ${\CC}$. This latter summand
is a (partial) contribution of the Eisenstein series $E_{a+2}$ of weight $a+2$.
We refer to the paper by Deligne \cite{Deligne}.

The relation just given is just one aspect of a deeper motivic relation; 
this aspect deals with the complex moduli space 
${\mathcal A}_1 \otimes {\CC}$; if we look
at ${\mathcal A}_1\otimes {\FF}_p$ we see another aspect.
For $\ell \neq p$ we have an isomorphism 
$$
H^i_c({\mathcal A}_1\otimes \overline{\FF}_p, {\VV}_a) {\buildrel \cong \over
\longrightarrow} 
H^i_c({\mathcal A}_1\otimes \overline{\QQ}_p, {\VV}_a)
$$
of ${\rm Gal}(\overline{\QQ}_p/{\QQ}_p)$-representations,
which bridges the gap between characteristic $0$ and characteristic $p$.
We can use this to see that for compactly supported 
\'etale $\ell$-adic cohomology with $\ell$ different from $p$, 
the trace of Frobenius on 
$H_c^1({\mathcal A}_1\otimes \overline{\FF}_p, {\VV}_a)$
equals $1$ plus 
the trace of the Hecke operator $T(p)$ on $S_{a+2}({\rm SL}(2,{\ZZ}))$, 
and this explains the identity (2). 
A more sophisticated version is that
$$
[H^1_c({\mathcal A}_1\otimes {\QQ}, {\VV}_a)]= S[a+2]+1 \, ,
\eqno(4)
$$
where the left hand side is viewed as a Chow motive with rational coefficients
and $S[k]$ denote the motive
associated by Scholl to the space of cusp forms of even weight $k>2$ on
${\rm SL}(2,{\ZZ})$. This incorporates both the Hodge theoretic 
and the Galois theoretic version.

But for elliptic curves and modular forms on ${\rm SL}(2,{\ZZ})$ we
have explicit knowledge and this way of mining 
information about modular forms by counting over finite fields might seem 
superfluous. Nevertheless, it is a practical method. Once one has a list
of all elliptic curves defined over ${\FF}_q$ up to isomorphism over ${\FF}_q$,
together with their number of points over ${\FF}_q$ and the order of their
automorphism groups, then one can easily calculate 
the trace of the Hecke operator $T(q)$
on the space $S_{k}({\rm SL}(2,{\ZZ}))$ for {\sl all} even weights $k>2$.

The situation changes drastically
if one considers curves of higher genus or abelian varieties of higher 
dimension and modular forms of higher degree. 
There our knowledge of modular forms is rather restricted and counting curves
over finite fields provides us with a lot of useful information
that is difficult to access otherwise.

We end this section with giving the relation between the cohomology 
of the local systems ${\VV}_a$ on ${\mathcal A}_1$ and the cohomology
of ${\mathcal M}_{1,n}$. 
The following beautiful formula due to Getzler
\cite{Getzler99} expresses the Euler characteristic 
$e_c({\mathcal M}_{1,n+1})$ in terms of the Euler characteristics
of the local systems ${\VV}_a$ in a concise way as a residue for $x=0$
in a formal expansion as follows:
$$
\frac{e_c({\mathcal M}_{1,n+1})}{n!}=
{\rm res}_0 \left[{L-x-L/x \choose n}
\sum_{k=1}^{\infty} \left( \frac{S[2k+2]+1}{L^{2k+1}})x^{2k}-1\right)
\cdot
(x-L/x) d x  \right]
$$
\end{section}
\begin{section}{Counting Curves of Genus Two}
The notion of elliptic curve allows two obvious generalizations:
one is that of an abelian variety of dimension $g>1$ 
and the other one is that of a curve of 
genus $g>1$. For $g=2$ these two generalizations are rather close.
The moduli space ${\mathcal M}_2$ of curves of genus $2$ admits an 
embedding in the moduli space ${\mathcal A}_2$ of principally polarized
abelian surfaces by the Torelli map, which
associates to a curve of genus $2$ its Jacobian. 
The image is an open part, the complement of the locus
${\mathcal A}_{1,1}$ of products of elliptic curves. 
The moduli spaces ${\mathcal M}_2$ and  ${\mathcal A}_{2}$ are 
defined over ${\ZZ}$.

The Hasse-Weil theorem tells us that for a curve $C$ of genus $2$ defined
over a finite field ${\FF}_q$ the action of Frobenius on 
$H^1_{\rm et}(C\otimes \overline{\FF}_q,{\QQ}_{\ell})$, with $\ell$ a prime different 
from the characteristic, is semi-simple and the eigenvalues $\alpha$
satisfy $\alpha \bar{\alpha}=q$.

The analogues of the notions that appeared in the preceding section 
are available.
We have the universal curve of genus $2$ over ${\mathcal M}_2$,
denoted by $\gamma: {\mathcal C}_2 \to
{\mathcal M}_2$, and
the universal principally abelian surface 
$\pi: {\mathcal X}_2 \to {\mathcal A}_2$. This gives rise to a local system
${\VV}:=R^1\pi_* {\QQ}_{\ell}$ on ${\mathcal A}_2$. 
This is a local system of rank $4$ and the pull back of this system 
under the Torelli morphism coincides with $R^1\gamma_*{\QQ}_{\ell}$.
The fibre of this local system ${\VV}$ over a point $[X]$ with $X$ 
a principally polarized abelian variety, is $H^1(X,{\QQ}_{\ell})$ and
this is a ${\QQ}_{\ell}$-vector space of dimension $4$ and ${\VV}$ 
is provided with a non-degenerate symplectic pairing 
${\VV}\times {\VV} \to {\QQ}_{\ell}(-1)$
that comes from the Weil pairing.

Instead of just considering the symmetric powers ${\rm Sym}^a{\VV}$ of 
${\VV}$, as we did for $g=1$, we can make more
local systems now. To every irreducible finite-dimensional representation 
of ${\rm Sp}(4,{\QQ})$, say of highest weight $\lambda=(a,b)$ with $a\geq b$, 
we can associate
a local system ${\VV}_{\lambda}$ by applying a Schur functor to ${\VV}$.
For $\lambda=(a,0)$ we recover ${\rm Sym}^a({\VV})$, and for 
example, ${\VV}_{(1,1)}$
is a $5$-dimensional local system occurring in $\wedge^2 {\VV}$.
A weight $\lambda=(a,b)$ is called {\sl regular} if $a>b>0$.

We then look at the  Euler characteristic
$$
\sum_{i=0}^{6} (-1)^i[H^i_c({\mathcal A}_2\otimes \overline{\QQ},
{\VV}_{\lambda})]\, ,
$$
where we consider the cohomology groups either as Hodge structures over
the complex numbers if we deal with complex cohomology over ${\mathcal A}_2
\otimes {\CC}$,
or as $\ell$-adic Galois representations when we consider $\ell$-adic
\'etale cohomology over ${\mathcal A}_2\otimes \overline{\QQ}$,
and the brackets indicate that the
sum is taken in a Grothendieck group of the appropriate category 
(Hodge structures or Galois representations). The information on the
cohomology over ${\FF}_p$ for all $p$ together gives
 the whole information over ${\QQ}$.

On the other hand the notion of modular form also generalizes. 
The moduli space
${\mathcal A}_2({\CC})$ of principally polarized complex abelian varieties 
can be represented by a quotient
$$
{\rm Sp}(4,{\ZZ})\backslash \mathfrak{H}_2
$$
with 
$\mathfrak{H}_2=
\{ \tau \in {\rm Mat}(2\times 2, {\CC}): \tau^t=\tau, {\rm Im}(\tau)>0\}$, 
the Siegel upper half space of degree~$2$. The symplectic group
${\rm Sp}(4,{\ZZ})$ acts on $\mathfrak{H}_2$ in the usual way by 
$$
\tau \mapsto (a\tau +b)(c \tau+d)^{-1} \qquad
\text{
for $\left( \begin{matrix} a & b \\ c & d \\ \end{matrix}\right) 
\in {\rm Sp}(4,{\ZZ})$. }
$$
A holomorphic function $f: \mathfrak{H}_2 \to W$ with $W$ a finite-dimensional
complex vector space that underlies a representation 
$\rho$ of ${\rm GL}(2,{\CC})$, 
is called a Siegel modular form of weight $\rho$ if
$f$ satisfies
$$
f((a\tau +b)(c \tau+d)^{-1})= \rho(c\tau+d) f(\tau) 
\qquad \text{ for all $\tau \in \mathfrak{H}_2$ and
$\left( \begin{matrix} a & b \\ c & d \\ \end{matrix}\right) 
\in {\rm Sp}(4,{\ZZ})$. }
$$
If $\rho$ is the $1$-dimensional representation $\det^k$, then $f$ is a
scalar-valued function and is called a classical Siegel modular form of
weight $k$. The Siegel modular forms of a given weight $\rho$ form a
finite-dimensional vector space $M_{\rho}({\rm Sp}(4,{\ZZ}))$. 
It contains a subspace $S_{\rho}({\rm Sp}(4,{\ZZ}))$ of cusp forms
characterized by a growth condition.

Without loss of generality we may consider only irreducible 
representations $\rho$ of ${\rm GL}(2)$. Such a representation
is of the form ${\rm Sym}^j(W)\otimes \det(W)^k$ with $W$ 
the standard representation. Therefore we shall use the 
notation $S_{j,k}({\rm Sp}(4,{\ZZ}))$
instead of $S_{\rho}({\rm Sp}(4,{\ZZ}))$, and similarly 
$M_{j,k}({\rm Sp}(4,{\ZZ}))$ for $M_{\rho}({\rm Sp}(4,{\ZZ}))$. 
We know that $M_{j,k}({\rm Sp}(4,{\ZZ}))$
vanishes if $j$ is odd or negative and also if $k$ is negative.
For the graded algebra of classical Siegel modular forms
$$
M=\oplus_k M_{0,k}({\rm Sp}(4,{\ZZ}))
$$
generators are known by work of Igusa. For a few cases of low values of $j$
we know generators for the $M$-module $\oplus_k M_{j,k}({\rm Sp}(4,{\ZZ}))$.

One also has a commutative algebra of Hecke operators acting on the spaces 
$M_{\rho}({\rm Sp}(4,{\ZZ}))$ and $S_{\rho}({\rm Sp}(4,{\ZZ}))$.
But in general we know much less than for genus $1$.

In order to get information about Siegel modular forms by counting curves 
of genus $2$ over finite fields 
one needs an analogue of the formula (2) (or (4)).

For genus $1$ we considered only the cohomology group $H^1$. 
It is known by work of Faltings that for a  local
system ${\VV}_{\lambda}$ with regular weight the cohomology groups 
$H^i_c({\mathcal A}_2\otimes {\QQ},{\VV}_{\lambda})$
vanish unless $i=3=\dim{\mathcal A}_2$. If $W_{\lambda}$ is an irreducible
representation of ${\rm Sp}(4,{\QQ})$ of highest weight $\lambda$,
then the Weyl character formula expresses the trace of an element
of ${\rm Sp}(4,{\QQ})$ as a symmetric function $\sigma_{\lambda}$
of the roots of its characteristic polynomial.

Since we can describe curves of genus $2$ very explicitly,
 we therefore consider for a curve $C$ of genus $2$ over ${\FF}_q$ 
with eigenvalues $\alpha_1,\bar{\alpha}_1,\alpha_2,\bar{\alpha}_2$
of Frobenius, i.e. 
such that
$$
\# C({\FF}_{q^n})= q^n+1-\alpha_1^n-\bar{\alpha}_1^n-\alpha_2^n -\bar{\alpha}_2^n\, ,
$$
the expression
$$
\frac{\sigma_{\lambda}(\alpha_1,\bar{\alpha}_1,\alpha_2,\bar{\alpha}_2)}{
\# {\rm Aut}_{{\FF}_q}(C)}
$$
and sum this over all isomorphism classes of genus $2$ curves defined 
over ${\FF}_q$. Here we are using the fact that each ${\FF}_q$-isomorphism 
class of genus $2$ curves defined over ${\FF}_q$ contains a curve
defined over ${\FF}_q$. In this way we find the analogue of the sum
$\sigma_a(q)$ defined in the preceding section. 
This gives us a way to calculate
the trace of Frobenius on the Euler characteristic 
of the cohomology of the local system ${\VV}_{a,b}$
on ${\mathcal M}_2\otimes \overline{\FF}_p$.
Define the (motivic) Euler characteristic
$$e_c({\mathcal M}_2\otimes {\QQ},{\VV}_{a,b}):=
\sum_{i=0}^6 (-1)^i [H^i({\mathcal M}_2\otimes {\QQ}, {\VV}_{a,b})]
$$
and similarly $e_c({\mathcal A}_2\otimes {\QQ},{\VV}_{a,b})$,
where the interpretation (Hodge structures or Galois modules)
depends on whether one takes complex cohomology or $\ell$-adic \'etale cohomology. 
We can then calculate the trace of Frobenius on 
$e_c({\mathcal M}_2 \otimes \overline{\FF}_p,{\VV}_{a,b})$ and 
$e_c({\mathcal A}_2\otimes \overline{\FF}_p,{\VV}_{a,b})$ by counting curves of 
genus $2$ over ${\FF}_p$.  
The difference between the two 
$$
e_c({\mathcal A}_2\otimes {\QQ},{\VV}_{a,b})-
e_c({\mathcal M}_2\otimes {\QQ},{\VV}_{a,b})=
e_c({\mathcal A}_{1,1}\otimes {\QQ}, {\VV}_{a,b})
$$
is the contribution 
from abelian surfaces that are products of two elliptic curves.
Or phrased differently, from stable curves of genus $2$ that consists of two
elliptic curves meeting in one point.  

How does this relate to the trace of Hecke operators on a space
$S_{\rho}({\rm Sp}(4,{\ZZ}))$? There is an analogue of the relation
(2), but the analogue of the term $1$ there is more complicated. 
Based on extensive calculations, in joint work with Carel Faber 
\cite{F-vdG} we 
formulated a conjecture
that is a precise analogue of (2).
We gave a formula for the Euler characteristic of the local system
${\VV}_{a,b}$ in the Grothendieck 
group of $\ell$-adic Galois representations.

The formula says that
$$
{\rm Trace}(T(p),S_{a-b,b+3}({\rm Sp}(4,{\ZZ})))= 
-{\rm Trace}(F_p,e_c({\mathcal A}_2\otimes \overline{\FF}_p,{\VV}_{a,b}))
+ {\rm Trace}(F_p,e_{2,\rm extra}(a,b))
$$
with $e_{2,\rm extra}(a,b)$ a correction term given by
$$
s_{a-b+2}+s_{a_b+4}(S[a-b+2]+1) {\mathbf L}^{b+1} +
\begin{cases} S[b+2]+1 & \text{$a$ even}\\
-S[a+3] & \text{$a$ odd.} \\
\end{cases}
$$
Here $s_k=\dim S_{k}({\rm SL}(2,{\ZZ}))$ and 
${\mathbf L}=h^2({\PP}^1)$ is the Lefschetz motive. 
The trace of Frobenius on ${\mathbf L}^k$ is $p^k$.

The conjecture has been proved by work of Weissauer for the regular case
and was completed by Petersen, see 
\cite{Weissauer1, Weissauer2, Petersen}.
One consequence is that the cohomology of the moduli spaces 
$\overline{\mathcal M}_{2,n}$ 
of stable $n$-pointed curves of genus $2$ is now completely known. 
It has also led to progress on the tautological rings of the moduli spaces 
${\mathcal M}_{2,n}$ by Petersen and Tommasi \cite{P-T}. 

This result allows us to calculate the traces of the Hecke operators
on spaces of classical and vector-valued Siegel modular forms.
The strategy to do this is by making a list of all Weil polynomials, 
that is, characteristic polynomials of Frobenius 
together with the frequency with which they occur 
if we go through all isomorphism classes, that is, if we run over
 ${\mathcal A}_2$. 
Once one has this list for a field ${\FF}_q$, one can compute 
the trace of the Hecke operator on the space of cusp forms 
$S_{j,k}({\rm Sp}(4,{\ZZ}))$ for {\sl all} pairs $(j,k)$ with $k\geq 3$.
We illustrate this with a few examples.

\begin{example}
The space $S_{0,35}({\rm Sp}(4,{\ZZ}))$ has dimension $1$ and is generated by the 
scalar-valued form $\chi_{35}$.
It corresponds to the case $(a,b)=(32,32)$. The eigenvalues of the Hecke
operators for $q\leq 37$ (and $q\neq 8,16,27, 32$) are given below. Note that
for $q$ a square the value differs from the usual one, see Definition 10.1
in \cite{B-F-vdG}. (The values for $q=p^r$ with $r\geq 3$  follow from those for
$q=p$ and $q=p^2$.)

\begin{footnotesize}
\smallskip
\vbox{
\bigskip\centerline{\def\quad{\hskip 0.6em\relax}
\def\quod{\hskip 0.5em\relax }
\vbox{\offinterlineskip
\hrule
\halign{&\vrule#&\strut\quod\hfil#\quad\cr
height2pt&\omit&&\omit&\cr
&$q$ && eigenvalue &\cr
\noalign{\hrule}
& $2$ && $-25073418240$ & \cr
& $3$ && $-11824551571578840$ & \cr
& $4$ && $ -203922016925674110976$ & \cr
& $5$ && $9470081642319930937500$ & \cr
& $7$ && $ -10370198954152041951342796400$ & \cr
& $9$ && $-270550647008022226363694871019974$ & \cr
& $11$ && $ -8015071689632034858364818146947656$ & \cr
& $13$ && $-20232136256107650938383898249808243380$ & \cr
& $17$ && $118646313906984767985086867381297558266980$ & \cr
& $19$ && $2995917272706383250746754589685425572441160$ & \cr
& $23$ && $ -1911372622140780013372223127008015060349898320$ & \cr
& $25$ && $-86593979298858393096680290648986986047363281250$ & \cr
& $29$ && $-2129327273873011547769345916418120573221438085460$ & \cr
& $31$ && $ -157348598498218445521620827876569519644874180822976$ & \cr
& $37$ && $-47788585641545948035267859493926208327050656971703460$ & \cr
} \hrule}
}}
\end{footnotesize}
\end{example}

\begin{example} Since explicitly known eigenvalues of Hecke operators
on Siegel modular forms are rather scarce, even for
scalar-valued forms of degree $2$, we give another example that shows how 
effective curve counting is. The space $S_{0,43}({\rm Sp}(4,{\ZZ}))$ is of 
dimension $1$ and generated
by a form $\chi_{43}= E_4^2 \chi_{35}$. We list the Hecke eigenvalues (with the same conventions
for prime powers as in the preceding one).

\begin{footnotesize}
\smallskip
\vbox{
\bigskip\centerline{\def\quad{\hskip 0.6em\relax}
\def\quod{\hskip 0.5em\relax }
\vbox{\offinterlineskip
\hrule
\halign{&\vrule#&\strut\quod\hfil#\quad\cr
height2pt&\omit&&\omit&\cr
&$q$ && eigenvalue &\cr
\noalign{\hrule}
& 2 && -4069732515840 & \cr
& 3 && -65782425978552959640 & \cr 
& 4 && -20941743921027137625128960 & \cr
& 5 && -44890110453445302863489062500 & \cr 
& 7 && -19869584791339339681013202023932400 & \cr
& 9 && -7541528134863704740446843276725979791820 & \cr
& 11 && 4257219659352273691494938669974303429235064 & \cr
& 13 && 1189605571437888391664528208235356059600166220 & \cr
& 17 && -1392996132438667398495024262137449361275278473925020 & \cr
& 19 && -155890765104968381621459579332178224814423111191589240 & \cr
& 23 && -128837520803382146891405898440571781609554910722934311120 & \cr
& 25 && 7099903749386561314439988230597055761986231311645507812500 & \cr
& 29 && 4716850092556381736632805755807948058560176106387507397101740 & \cr
& 31 && 3518591320768311083473550005851115474157237215091087497259584 & \cr
& 37 && -80912457441638062043356244171113052936003605371913289553380964260 & \cr
} \hrule}
}}
\end{footnotesize}
\end{example}

\bigskip
\begin{example}
The first cases where one finds a vector-valued cusp form that is not a lift
from elliptic modular forms
are the cases $(j,k)=(6,8)$ and $(4,10)$. We give the eigenvalues.
We also give the eigenvalues for $(j,k)=(34,4)$. In all these cases the
space of cusp forms is $1$-dimensional.

\begin{footnotesize}
\smallskip
\vbox{
\bigskip\centerline{\def\quad{\hskip 0.6em\relax}
\def\quod{\hskip 0.5em\relax }
\vbox{\offinterlineskip
\hrule
\halign{&\vrule#&\strut\quod\hfil#\quad\cr
height2pt&\omit&&\omit&&\omit&&\omit&\cr
&$q$ && $(6,8)$ && $(4,10)$&& $(34,4)$ &\cr
\noalign{\hrule}
& $2$ && $0$ && $-1680$ && $-633600$ &\cr
& $3$ && $-27000$ && $55080$ && $91211400$ & \cr
& $4$ && $409600$ && $-700160$ && $271415050240$ & \cr
& $5$ && $2843100$ &&$-7338900$ && $11926488728700$ & \cr
& $7$ && $-107822000$ && $609422800$ && $6019524504994000$ & \cr
& $9$ && $333371700$ && $1854007380$ && $-1653726849656615820$ & \cr
& $11$ && $3760397784$ && $25358200824$ && $-121499350185684258216$ & \cr
& $13$ && $9952079500$ && $-263384451140$ && $655037831218999528300$ & \cr
& $17$ && $243132070500$ && $-2146704955740$ && $714735598649071209833700$ & \cr
& $19$ && $595569231400$ && $43021727413960$ && $-3644388446450362098497240$ & \cr
& $23$  && $-6848349930000$ && $-233610984201360$ && $179014316167538651075065200$ & \cr
& $25$ && $-15923680827500$ && $-904546757727500$ && $-1338584707016863344819747500$ & \cr
& $29$ && $53451678149100$ && $-545371828324260$ && $52292335454052856173814993740$ & \cr
& $31$ && $234734887975744$ && $830680103136064$  && $ -256361532431714633455270321856$ & \cr
& $37$ && $448712646713500$ && $11555498201265580$  && $-826211657019923608686387368900$ &  \cr
} \hrule}
}}
\end{footnotesize}

\end{example} 

The fact that one can calculate these eigenvalues has motivated Harder
to make an idea about congruences between elliptic modular forms and Siegel
modular forms of degree $2$ concrete and formulate a conjecture about such
congruences. Already many years ago, 
Harder had the idea that there should be congruences between
the Hecke eigenvalues of elliptic modular forms and Siegel modular forms 
of genus $2$ modulo a prime that divides a critical value of the L-function
of the elliptic modular form, 
but the fact that the genus $2$ eigenvalues
could be calculated spurred him to make his ideas more concrete. 
He formulated his conjecture in \cite{Harder1}.
 
If $f=\sum a(n) q^n \in S_k({\rm SL}(2,{\ZZ}))$ is a normalized ($a(1)=1$) 
elliptic modular cusp form with 
L-function $\sum a(n)n^{-s}$ and $\Lambda(f,s)=
(\Gamma(s)/(2\pi)^s) L(f,s)$ that satisfies
the functional equation $\Lambda(f,s)=(-1)^{k/2} \Lambda(f,k-s)$, then the values
$\Lambda(f,r)$ with $k/2 \leq r \leq k-1$ are called the {\sl critical values}.
According to Manin and Vishik there are real numbers $\omega_{\pm}(f)$ 
with the property that all values $\Lambda^{\prime}(f,r):=
\Lambda(f,r)/\omega_{+}(f)$ for 
$r$ even (resp.\ $\Lambda^{\prime}(f,r):=\Lambda(f,r)/\omega_{-}(f)$ for $r$ odd) lie in 
${\QQ}_f={\QQ}(a(n): n\in {\ZZ}_{\geq 1})$, the field of eigenvalues 
$\lambda_p(f)=a(p)$ of the Hecke operators.
If $\ell$ is a prime in ${\QQ}_f$ lying above $p$ it is called ordinary
if $a(p)\not\equiv 0 \, (\bmod\,  \ell)$.

Harder's conjecture says the following.

\begin{conjecture} {\rm (Harder \cite{Harder1})} Let $a>b$ be natural numbers and $f\in S_{a+b+4}({\rm SL}(2,{\ZZ}))$ be an eigenform. If $\ell$ is an ordinary prime in the field
${\QQ}_f$ of Hecke eigenvalues $\lambda_p(f)$ of $f$ 
and $\ell^s$ with $s\geq 1$ divides the critical  value $\Lambda^{\prime}(f,a+3)$,
 then there exists an eigenform $F \in S_{a-b,b+3}({\rm Sp}(4,{\ZZ}))$ with 
Hecke eigenvalues $\lambda_p(F)$ satisfying
$$
\lambda_p(F) \equiv p^{a+2}+\lambda_p(f)+p^{b+1} \, (\bmod \ell^s)
$$
in the ring of integers of the compositum of the fields ${\QQ}_f$ and ${\QQ}_F$ 
of Hecke eigenvalues of $f$ and $F$ for all primes $p$.
\end{conjecture}

The counting of curves
over finite fields provided a lot of evidence for his conjecture. We 
give one example.

Let $(a,b)=(20,4)$ and let $f\in S_{28}({\rm SL}(2,{\ZZ}))$ be the normalized
eigenform. This form has eigenvalues in the field ${\QQ}(\sqrt{d})$ 
with $d=18209$.
We have $f=\sum_{n\geq 1} a(n) \, q^n$ with 
$$
f= q+(-4140-108\, \sqrt{d})q^2 +(-643140-20737\, \sqrt{d})q^3 + \cdots
$$
with 
$$
a(37)= \lambda_{37}(f)= 933848602341412283390+4195594851869555712\,\sqrt{d} .
$$
The critical value of $\Lambda(f,22)$ is divisible by the ordinary prime $367$.
Harder's conjecture claims that there should be a congruence.
Indeed, the space $S_{16,7}({\rm Sp}(4,{\ZZ}))$ has dimension $1$
and is thus generated by a Hecke eigenform $F$
and our results give  the eigenvalue 
$$
\lambda_{37}(F)=-1845192652253792587940.
$$ 
The prime 367 splits in ${\QQ}(\sqrt{18209})$ as $367=\pi\cdot \pi'$
with $\pi =(367,260+44 \, \sqrt{d})$.
The reader may check that indeed we 
have the congruence 
$$
\lambda_{37}(F)\equiv 37^{22}+ a(37)+37^5 \, (\bmod \pi) \, .
$$
Thus the counting of curves provides evidence for these conjectures.
For more details see \cite{vdG1, Harder1, B-F-vdG}.
\end{section}

\begin{section}{Counting Curves of Genus Three}
Like for genus $2$, the moduli spaces ${\mathcal M}_3$ of curves of genus $3$
and ${\mathcal A}_3$ of principally polarized abelian varieties of dimension
$3$ are rather close; in this case the Torelli map is a morphism 
$t: {\mathcal M}_3 \to {\mathcal A}_3$ of 
Deligne-Mumford stacks of degree $2$. This is due to the fact that
every abelian variety $X$ has an automorphism of order $2$ given by $-1_X$,
while the generic curve of genus $3$ has no non-trivial automorphisms.
The universal families $\pi: {\mathcal X}_3 \to {\mathcal A}_3$ and
$\gamma: {\mathcal C}_3 \to {\mathcal M}_3$ define local systems
${\VV}:=R^1\pi_* {\QQ}_{\ell}$ and $R^1\gamma_* {\QQ}_{\ell}$ with
the pull back $t^* {\VV}=R^1\gamma_*{\QQ}_{\ell}$. 
The local system ${\VV}$ carries a non-degenerate symplectic pairing
${\VV} \times {\VV} \to {\QQ}_{\ell}(-1)$ and again we find for each
irreducible representation of ${\rm Sp}(6,{\QQ})$ of highest weight $\lambda=
(a,b,c)$ with $a\geq b \geq c\geq 0$ a local system ${\VV}_{\lambda}$.
We are interested
in 
$$
e_c({\mathcal A}_3,{\VV}_{\lambda})=\sum_{i=0}^{12} (-1)^i [H^i_c({\mathcal A}_3,{\VV}_{\lambda})]\, ,
$$
again viewed in a Grothendieck group of Hodge structures or Galois representations.

What does the trace of Frobenius on this Euler characteristic tell us about traces of
Hecke operators on Siegel modular forms? Here a Siegel modular form is a 
holomorphic function $f: \mathfrak{H}_3 \to W$ with $W$ a finite-dimensional 
complex representation $\rho$ of ${\rm GL}(3,{\CC})$ satisfying
$$
f((a\tau +b)(c \tau+d)^{-1})= \rho(c\tau+d) f(\tau) 
\qquad \text{ for all $\tau \in \mathfrak{H}_3$ and
$\left( \begin{matrix} a & b \\ c & d \\ \end{matrix}\right) 
\in {\rm Sp}(6,{\ZZ})$ }.
$$
If $\rho$ is an irreducible representation of ${\rm GL}(3)$ of highest
weight $(\alpha,\beta,\gamma)$ with $\alpha\geq \beta \geq \gamma$,
then the corresponding space of modular forms (resp.\ cusp forms)
is denoted by $M_{i,j,k}({\rm Sp}(6,{\ZZ}))$ (resp.\ by
$S_{i,j,k}({\rm Sp}(6,{\ZZ}))$) 
and their weight (in the sense of modular forms)
is denoted with $(i,j,k)=(\alpha-\beta,\beta-\gamma,\gamma)$.

In joint work with Bergstr\"om and Faber \cite{B-F-vdG} we formulated 
a conjecture relating the trace of the
Hecke operator on a space of vector-valued Siegel modular forms
with the counts of curves. It was based on extensive calculations
using counting of curves.
It says
\begin{conjecture}\label{conjg3} 
For $\lambda=(a,b,c)$ the trace of the Hecke operator $T(p)$ on 
the space of cusp forms 
$S_{a-b,b-c,c+4}({\rm Sp}(6,{\ZZ}))$ is given by the trace of Frobenius
on $e_c({\mathcal A}_3,{\VV}_{\lambda})$ minus a correction term
$e_{3,{\rm extra}}(a,b,c)$ given by
$$
\begin{aligned}
e_{3,{\rm extra}} =
-e_c({\mathcal A}_2,{\VV}_{a+1,b+1}) + e_c({\mathcal A}_2,{\VV}_{a+1,c})
&-e_c({\mathcal A}_2,{\VV}_{b,c}) \\
-e_{2,{\rm extra}}(a+1,b+1) \otimes S[c+2]
&+e_{2,{\rm extra}}(a+1,c) \otimes S[b+3]\\
& \qquad -e_{2,{\rm extra}}(b,c) \otimes S[a+4]  \\
\end{aligned}
$$
\end{conjecture}

The evidence for this conjecture is overwhelming. It fits with all we know about classical Siegel modular forms. The dimensions fit with the numerical
Euler characteristics (replacing $[H^i_c({\mathcal A}_3,{\VV}_{\lambda})]$
by its dimension). Moreover, the answers that we find by counting turn out
to be integers, which is already quite a check, 
as we are summing rational numbers due to the factors $1/\# {\rm Aut}_{{\FF}_q}(C)$.
These results also fit with very recent (conjectural) results concerning Siegel
modular forms obtained by the Arthur trace formula, see \cite{Chenevier-Renard}.

In order to show that it leads to very concrete results we give an
illustration.

\begin{example}
The lowest weight examples of cusp forms that are not lifts occur in
weights $(3,3,7),(4,2,8)$ and $(2,6,6)$. In these cases
the spaces $S_{i,j,k}({\rm Sp}(6,{\ZZ}))$
are $1$-dimensional. We give the (conjectured) Hecke eigenvalues.

\begin{footnotesize}
\smallskip
\vbox{
\bigskip\centerline{\def\quad{\hskip 0.6em\relax}
\def\quod{\hskip 0.5em\relax }
\vbox{\offinterlineskip
\hrule
\halign{&\vrule#&\strut\quod\hfil#\quad\cr
height2pt&\omit&&\omit&&\omit&& \omit &\cr
&$p\backslash (i,j,k)$ && $(3,3,7)$ && $(4,2,8)$&& $(2,6,6)$&\cr
\noalign{\hrule}
& $2$ && $1080$ && $9504$ && $5184$ & \cr
& $3$ && $181440$ && $970272$  && $-127008$ & \cr
& $4$ && $15272000$ && $89719808$  && $62394368$ &  \cr
& $5$ && $368512200$ && $-106051896$  && $2126653704$ & \cr
& $7$ && $13934816000$ && $112911962240$  && $86958865280$ & \cr
& $8$ && $-15914672640$ && $1156260593664$  && $32296402944$ & \cr
& $9$ && $483972165000$ && $5756589166536$  && $1143334399176$ & \cr
& $11$ && $424185778368$ && $44411629220640$  && $64557538863840$ & \cr
& $13$ && $26955386811080$ && $209295820896008$  && $-34612287925432$ & \cr
& $16$ && $1224750814466048$ && $-369164249202688$  && $12679392014630912$ & \cr
& $17$ && $282230918895240$ && $1230942201878664$  && $7135071722206344$ & \cr
& $19$ && $5454874779704000$ && $51084504993278240$ && $46798706961571040$ & \cr
} \hrule}
}}
\end{footnotesize}

\end{example}
\begin{example} As stated, our data allow the calculation for the trace of the
Hecke operator $T(q)$ for $q\leq 19$ for all weights if Conjecture \ref{conjg3}
is granted. 
Here we present the case of weight $(60,0,4)$. The space 
$S_{60,0,4}({\rm Sp}(6,{\ZZ}))$ is of dimension $1$ and a generator 
is not a lift. We give the 
Hecke eigenvalues, where
for prime powers we use the convention of 10.1 in \cite{B-F-vdG}.

\begin{footnotesize}
\smallskip
\vbox{
\bigskip\centerline{\def\quad{\hskip 0.6em\relax}
\def\quod{\hskip 0.5em\relax }
\vbox{\offinterlineskip
\hrule
\halign{&\vrule#&\strut\quod\hfil#\quad\cr
height2pt&\omit&& \omit &\cr
&$p$ && {\rm eigenvalue} &\cr
\noalign{\hrule}
& $2$ && $1478987712$ & \cr
& $3$ && $-2901104577414432$ & \cr
& $4$ && $-81213310977988096000$ & \cr
& $5$ && $17865070879279088017800$ & \cr 
& $7$ && $-6212801311610929434173542528$ & \cr 
& $8$ && $-1127655095344679889821203955712$ & \cr
& $9$ && $5614158763137782860896126573000$ &\cr
& $11$ && $-1849697485178583502997203666501152$ & \cr
& $13$ && $2477960171489248682447718208861099208$ &\cr
& $16$ && $-8941917317486628689603624398015726354432$ & \cr 
& $17$ && $-73908079488243072323266509093278640761208$ & \cr
& $19$ && $592331726239601530766675208936630486956000$ & \cr
} \hrule}
}}
\end{footnotesize}
\end{example}

As for genus $2$, 
the heuristics of counting has led to new conjectured liftings 
to modular forms of genus $3$, 
to new Harder type congruences and other results. 
We refer to \cite{B-F-vdG}. 
\end{section}

\begin{section}{Other Cases}
For genus $g\geq 4$ the dimension of the moduli space ${\mathcal A}_g$
of principally polarized abelian varieties of dimension $g$ is larger than the
dimension of the moduli space ${\mathcal M}_g$ of curves of genus $g$.
This means that one cannot use the Torelli map $t: {\mathcal M_g \otimes {\FF}_q
} \to {\mathcal A}_g\otimes {\FF}_q$, which associates to a curve its Jacobian
variety,
to enumerate all principally polarized abelian varieties of dimension $g$
over ${\FF}_q$. For genus $g=4$ or $5$ one might consider instead the 
Prym varieties of double \'etale covers of curves of genus $g+1$, 
but enumerating
these double covers is already considerably more difficult. And for $g\geq 7$
the moduli spaces ${\mathcal A}_g$ are of general type, hence not unirational
and therefore cannot be parametrized by open parts of affine or projective
spaces. Nevertheless, there are other families of curves and abelian varieties
to which the method of counting over finite fields can be applied.

In \cite{Shimura1} Shimura describes a number of moduli spaces that over the
complex numbers have a complex ball as universal cover and are rational
varieties (birationally equivalent to projective space). In all these
cases these are moduli spaces of curves that are described as covers of the 
projective line. One such case concerns triple Galois covers of genus $3$ 
of the projective line. If the characteristic of the field is not $3$, then
such a curve can be given as $y^3=f(x)$ with $f\in k[x]$ a degree $4$
polynomial with distinct zeros. The Jacobians of such curves
are abelian threefolds with multiplication by $F={\QQ}(\sqrt{-3})$
induced by the action of the Galois automorphism $\alpha$ of the curve
of order $3$.
The moduli of such abelian threefolds over ${\CC}$ are described by
an arithmetic quotient of the complex $2$-ball by a discrete subgroup
of the algebraic group of similitudes
$G=\{ g \in {\rm GL}(3,F): h(gz,gz)=\eta(g)h(z,z)\}$
of a non-degenerate hermitian form 
$h=z_1\bar{z}_2+z_2\bar{z}_1+z_3\bar{z}_3$ 
on $F^3$, 
where the bar refers to the Galois automorphism of $F$. 
In fact, the discrete subgroup is the group $\Gamma[\sqrt{-3}]$
$$
\{ g \in {\rm GL}(3,O_F): h(gz,gz)=h(z,z), g\equiv 1 (\bmod \sqrt{-3}) \}
$$

On our moduli space ${\mathcal M}$ defined over the ring of integers $O_F[1/3]$
of $F$ with $3$ inverted we have a universal family $\pi: {\mathcal C} \to {\mathcal M}$ and hence we get again a local system ${\VV}=R^1\pi_* {\QQ}$
or $R^1\pi_* {\QQ}_{\ell}$. This is a local system of rank $6$ provided with
a non-degenerate alternating pairing ${\VV}\times {\VV} \to {\QQ}(-1)$.
The action of $\alpha$ on the cohomology gives rise to a splitting of 
${\VV}$ as a direct sum of two local systems of rank $3$ over $F$:
${\VV}\otimes F={\WW}\oplus {\WW}'$. The non-degenerate pairing implies that
${\WW}'\cong {\WW}^{\vee} \otimes F(-1)$, where we denote by ${\WW}^{\vee}$ the $F$-linear dual.
From these basic local systems ${\WW}$, ${\WW}'$ 
we can obtain for each irreducible
representation $\rho$ of ${\rm GL}(3)$ local systems that appear as the analogues
of the local systems ${\VV}_a$ for $g=1$ and ${\VV}_{\lambda}$ for $g=2$ and $3$.

The role of the Siegel modular forms is now taken by so-called 
Picard modular forms. In fact, identifying $G({\QQ})$ with the matrix subgroup
of ${\rm GL}(3,F)$ this group acts on the domain $B=
\{ (u,v) \in {\CC}^2: 2 {\rm Re}(v)+|u|^2 <0 \}$ (isomorphic to a complex ball)
by 
$$
(u,v) \mapsto \left( \frac{g_{31}v+g_{32}+g_{33}u}{g_{21}v+g_{22}+g_{23}u},
\frac{g_{11}v +g_{12}+g_{13}u}{g_{21}v+g_{22}+g_{23}u} \right) .
$$
For $g=(g_{ij})\in G$ we let
$$
j_1(g,u,v)=g_{21}v+g_{22}+g_{23}u
$$
and
$$
j_2(g,u,v^{-1})=\det(g)^{-1} \left( 
\begin{matrix}
G_{32}u+G_{33} & G_{12}u+G_{13}\\
G_{12} u +G_{13} & G_{12}v+ G_{11} \\ 
\end{matrix} \right)
$$
with $G_{ij}$ the minor of $g_{ij}$. Then a (vector-valued) 
Picard modular form
of weight $(j,k)$ on our discrete subgroup $\Gamma[\sqrt{-3}]$ is a 
holomorphic map $f: B \to {\rm Sym}^j({\CC}^2)$ satisfying 
$$
f(g\cdot (u,v))=j_1(g,u,v)^k {\rm Sym}^j(j_2(g,u,v)) f(u,v)
$$
for all $g \in \Gamma[\sqrt{-3}]$.

In joint work with Bergstr\"om we analyzed the  Euler characteristic of
compactly supported cohomology
of local systems in this case by extensive counting over finite fields 
and came forward with conjectures that describe
the Euler characteristics 
of these local systems and  the traces of Hecke operators
on the corresponding spaces of Picard modular forms, see \cite{B-vdG}.
These conjectures guided work of Cl\'ery and van der Geer to construct
the vector-valued modular forms and to find generators for modules of 
such vector-valued Picard modular forms. We refer to \cite{C-vdG1}.

\bigskip
One of the charms of the subject of curves over finite fields is, that it is 
relatively easily accessible without requiring sophisticated techniques
and amenable to direct calculations. Although it arose late, it is
intimately connected to very diverse array of subdisciplines of mathematics.
I hope to have convinced the reader that it is also a
 wonderful playground to find heuristically new phenomena and
patterns that can help other areas of mathematics. 

\smallskip
{\sl Acknowledgement} The author thanks Jonas Bergstr\"om, Fabien Cl\'ery 
and the referees for some helpful remarks.
\end{section}

\begin{section}{Tables}

The following two tables summarize the status quo as contained in the tables
of the website {\tt www.manypoints.org} for the function $N_q(g)$ 
for $1 \leq g \leq 50$ and $q$ equal to a small power of $2$ or $3$.
It gives either one value for $N_q(g)$, or an interval $[a,b]$ 
(denoted as $a-b$ in the tables) in which
$N_q(g)$ is supposed to lie, or an entry $-b$ if $b$ is the best upper 
bound known for $N_q(g)$ and no curve with at least $[b/\sqrt{2}]$ 
rational points is known, see \cite{vdG-vdV3}. 
\vfill\eject
%
\font\tablefont=cmr8
\def\quad{\hskip 0.6em\relax}
\def\quod{\hskip 0.6em\relax}
\def\vhop{
    height2pt&\omit&&\omit&&\omit&&\omit&&\omit&&\omit&&\omit&&\omit&\cr}
\noindent{\bf Table p=2.}
$$
\vcenter{
\tablefont
\lineskip=1pt
\baselineskip=10pt
\lineskiplimit=0pt
\setbox\strutbox=\hbox{\vrule height .7\baselineskip
                                depth .3\baselineskip width0pt}%
\offinterlineskip
\hrule
\halign{&\vrule#&\strut\quod\hfil#\quad\cr
\vhop
&$g\backslash q$&&2&&4&&8&&16&&32&&64&&128&\cr
\vhop
\noalign{\hrule}
\vhop
&1&&5&&9&&14&&25&&44&&81&&150&\cr
&2&&6&&10&&18&&33&&53&&97&&172&\cr
&3&&7&&14&&24&&38&&64&&113&&192&\cr
&4&&8&&15&&25&&45&&71--72&&129&&215&\cr
&5&&9&&17&&29&&49--53&&83--85&&140--145&&227--234&\cr
&6&&10&&20&&33--34&&65&&86--96&&161&&243--256&\cr
&7&&10&&21&&34--38&&63--69&&98--107&&177&&262--283&\cr
&8&&11&&21--24&&35--42&&63--75&&97--118&&169--193&&276--302&\cr
&9&&12&&26&&45&&72--81&&108--128&&209&&288--322&\cr
&10&&13&&27&&42--49&&81--86&&113--139&&225&&296--345&\cr
\noalign{\hrule}
&11&&14&&26--29&&48--53&&80--91&&120--150&&201--235&&294--365&\cr
&12&&14--15&&29--31&&49--57&&88--97&&129--160&&257&&321--388&\cr
&13&&15&&33&&56--61&&97--102&&129--171&&225--267&&--408&\cr
&14&&16&&32--35&&65&&97--107&&146--182&&257--283&&353--437&\cr
&15&&17&&35--37&&57--67&&98--112&&158--193&&258--299&&386--454&\cr
&16&&17--18&&36--38&&56--70&&95--118&&147--204&&267--315&&--476&\cr
&17&&18&&40&&63--73&&112--123&&154--211&&--331&&--499&\cr
&18&&18--19&&41--42&&65--77&&113--128&&161--219&&281--347&&--519&\cr
&19&&20&&37--43&&60--80&&129--133&&172--227&&315--363&&--542&\cr
&20&&19--21&&40--45&&76--83&&127--139&&177--235&&342--379&&--562&\cr
\noalign{\hrule}
&21&&21&&44--47&&72--86&&129--144&&185--243&&281--395&&--591&\cr
&22&&21--22&&42--48&&74--89&&129--149&&--251&&321--411&&--608&\cr
&23&&22--23&&45--50&&68--92&&126--155&&--259&&--427&&--630&\cr
&24&&23&&49--52&&81--95&&129--161&&225--266&&337--443&&513--653&\cr
&25&&24&&51--53&&86--97&&144--165&&--274&&408--459&&--673&\cr
&26&&24--25&&55&&82--100&&150--170&&--282&&425--475&&--696&\cr
&27&&24--25&&52--56&&96--103&&156--176&&213--290&&416--491&&--716&\cr
&28&&25--26&&54--58&&97--106&&145--181&&257--297&&513&&577--745&\cr
&29&&26--27&&52--60&&97--109&&161--186&&227--305&&--523&&--761&\cr
&30&&25--27&&53--61&&96--112&&162--191&&273--313&&464--535&&609--784&\cr
\noalign{\hrule}
&31&&27--28&&60--63&&89--115&&168--196&&--321&&450--547&&578--807&\cr
&32&&27--29&&57--65&&90--118&&--201&&--328&&--558&&--827&\cr
&33&&28--29&&65--66&&97--121&&193--207&&--336&&480--570&&--850&\cr
&34&&27--30&&65--68&&98--124&&183--212&&--344&&462--581&&--870&\cr
&35&&29--31&&64--69&&112--127&&187--217&&253--351&&510--593&&--899&\cr
&36&&30--31&&64--71&&112--130&&185--222&&--359&&490--604&&705--914&\cr
&37&&30--32&&66--72&&121--132&&208--227&&--367&&540--615&&--938&\cr
&38&&30--33&&64--74&&129--135&&193--233&&291--375&&518--627&&--961&\cr
&39&&33&&65--75&&120--138&&194--238&&--382&&494--638&&--981&\cr
&40&&32--34&&75--77&&103--141&&225--243&&293--390&&546--649&&--1004&\cr
\noalign{\hrule}
&41&&34--35&&72--78&&118--144&&220--248&&308--397&&560--661&&--1024&\cr
&42&&33--35&&75--80&&129--147&&209--253&&307--405&&574--672&&--1053&\cr
&43&&34--36&&72--81&&116--150&&226--259&&306--412&&546--683&&--1068&\cr
&44&&33--37&&68--83&&130--153&&226--264&&325--420&&516--695&&--1092&\cr
&45&&36--37&&80--84&&144--156&&242--268&&313--427&&572--706&&--1115&\cr
&46&&36--38&&81--86&&129--158&&243--273&&--435&&585--717&&--1135&\cr
&47&&36--38&&73--87&&126--161&&--277&&--443&&598--728&&--1158&\cr
&48&&35--39&&80--89&&128--164&&243--282&&--450&&564--739&&--1178&\cr
&49&&36--40&&81--90&&130--167&&240--286&&--458&&624--751&&913--1207&\cr
&50&&40&&91--92&&130--170&&255--291&&--465&&588--762&&--1222&\cr
\vhop
}
\hrule
}
$$
%
\vfill\eject

\font\tablefont=cmr8
\def\quad{\hskip 0.6em\relax}
\def\quod{\hskip 0.6em\relax}
\def\vhop{height2pt&\omit&&\omit&&\omit&&\omit&&\omit&\cr}
\noindent{\bf Table p=3.}
\bigskip
$$
\vcenter{
\tablefont
\lineskip=1pt
\baselineskip=10pt
\lineskiplimit=0pt
\setbox\strutbox=\hbox{\vrule height .7\baselineskip
                                depth .3\baselineskip width0pt}%
\offinterlineskip
\hrule
\halign{&\vrule#&\strut\quod\hfil#\quad\cr
height2pt&\omit&&\omit&&\omit&&\omit&&\omit&\cr
&$g\backslash q$&&3&&9&&27&&81&\cr
height2pt&\omit&&\omit&&\omit&&\omit&&\omit&\cr
\vhop
\noalign{\hrule}
\vhop
height2pt&\omit&&\omit&&\omit&&\omit&&\omit&\cr
&1&&7&&16&&38&&100&\cr
&2&&8&&20&&48&&118&\cr
&3&&10&&28&&56&&136&\cr
&4&&12&&30&&64&&154&\cr
&5&&13&&32--35&&72--75&&167--172&\cr
&6&&14&&35--38&&76--84&&190&\cr
&7&&16&&40--43&&82--95&&180--208&\cr
&8&&17--18&&42--46&&92--105&&226&\cr
&9&&19&&48--50&&99--113&&244&\cr
&10&&20--21&&54&&94--123&&226--262&\cr
\noalign{\hrule}
&11&&21--22&&55--58&&100--133&&220--280&\cr
&12&&22--23&&56--61&&109--143&&298&\cr
&13&&24--25&&64--65&&136--153&&298--312&\cr
&14&&24--26&&56--69&&--163&&278--330&\cr
&15&&28&&64--73&&136--170&&292--348&\cr
&16&&27--29&&74--76&&144--178&&370&\cr
&17&&28--30&&74--80&&128--184&&288--384&\cr
&18&&28--31&&68--84&&148--192&&306--400&\cr
&19&&32&&84--88&&145--199&&--418&\cr
&20&&30--34&&70--91&&--206&&--436&\cr
\noalign{\hrule}
&21&&32--35&&88--95&&163--213&&352--454&\cr
&22&&33--36&&78--98&&--220&&370--472&\cr
&23&&33--37&&92--101&&--227&&--490&\cr
&24&&31--38&&91--104&&208--234&&370--508&\cr
&25&&36--40&&96--108&&196--241&&392--526&\cr
&26&&36--41&&110--111&&200--248&&500--544&\cr
&27&&39--42&&104--114&&208--255&&--562&\cr
&28&&37--43&&105--117&&--262&&--580&\cr
&29&&42--44&&104--120&&196--269&&--598&\cr
&30&&38--46&&91--123&&196--276&&551--616&\cr
\noalign{\hrule}
&31&&40--47&&120--127&&--283&&460--634&\cr
&32&&40--48&&93--130&&--290&&--652&\cr
&33&&48--49&&128--133&&220--297&&576--670&\cr
&34&&46--50&&111--136&&--304&&594--688&\cr
&35&&47--51&&119--139&&--311&&612--706&\cr
&36&&48--52&&118--142&&244--318&&730&\cr
&37&&52--54&&126--145&&236--325&&648--742&\cr
&38&&--55&&111--149&&--332&&629--755&\cr
&39&&48--56&&140--152&&271--340&&730--768&\cr
&40&&56--57&&118--155&&273--346&&663--781&\cr
\noalign{\hrule}
&41&&50--58&&140--158&&--353&&680--794&\cr
&42&&52--59&&122--161&&280--360&&697--807&\cr
&43&&56--60&&147--164&&--367&&672--821&\cr
&44&&47--61&&119--167&&278--374&&--834&\cr
&45&&54--62&&136--170&&--380&&704--847&\cr
&46&&60--63&&162--173&&--387&&720--859&\cr
&47&&54--65&&154--177&&299--394&&690--872&\cr
&48&&55--66&&163--180&&325--401&&752--885&\cr
&49&&64--67&&168--183&&316--408&&768--898&\cr
&50&&63--68&&182--186&&312--415&&784--911&\cr
\vhop
}
\hrule
}
$$
\end{section}

\end{document}